\documentclass[12pt]{article}  
\pdfoutput=1 	
\usepackage{geometry}                		
\usepackage[parfill]{parskip}    		
\usepackage{graphicx}			
\usepackage[hyperindex=true,pagebackref,colorlinks=true,pdfpagemode=none,urlcolor=blue,linkcolor=blue,citecolor=blue,pdfstartview=FitH]{hyperref}
\usepackage{amsmath,amsfonts}
\usepackage{color}
\usepackage{amsthm}
\usepackage{boolexpr}
\usepackage{amssymb,latexsym,pstricks,pst-plot}
\usepackage[english]{babel}
\usepackage{pgf}
\usepackage{tikz}
\usetikzlibrary{arrows,automata}
\usepackage[UKenglish]{isodate}
\usepackage[T1]{fontenc}
\usepackage{listings}
\usepackage{xcolor}
\usepackage{bm}
\usepackage{mathabx}
\usepackage{enumitem}
\usepackage{mathtools}
\allowdisplaybreaks

\usepackage{amssymb}

\theoremstyle{plain}

\newtheorem{theorem}{Theorem}

\theoremstyle{definition}

\theoremstyle{remark}

\title{Chung-Graham and Zeckendorf representations}
\author{Rob Burns}

\begin{document}
\maketitle
\begin{abstract}
We examine the relationship between the Chung-Graham and Zeckendorf representations of an integer using the software package {\tt Walnut}.
\end{abstract}

\section{Introduction}
\label{intro}

There are numerous ways of representing positive integers as a linear sum of a fixed set of other positive integers. The most familiar representation involves writing an integer $n$ in what is called base $b$:
$$
n = \sum_i a_i b^i \,\,\, \text{where } \,\, 0 \leq a_i < b.
$$
Examples include the binary and decimal number systems.  

We denote the Fibonacci numbers by $F_i, \,\, i\geq 0$. They are defined by the recurrence
$$
F_0 = 0, \,\, F_1 = 1, \,\, F_{n+1} = F_{n} + F_{n-1} \,\, \text{ for } \,\, n \geq 1.
$$
The Fibonacci numbers can be used as a basis for representing integers. Some examples are provided in the paper \cite{Shallit:2023aa} by Shallit and Shan and in the paper by Gilson et al. \cite{Gilson:2020aa}. Shan's thesis describes how to use the pre-installed representations in {\tt Walnut} to prove that other regular numeration systems are complete and unambiguous \cite{shan2024}. \textit{Regular} here means that the system is recognised by a finite automaton. \textit{Complete and unambiguous} means the system has exactly one way of representing every positive integer. The most well known system involving the Fibonacci numbers is the Zeckendorf representation, which was published by Lekkerkerker in 1952 \cite{lekk} and by Zeckendorf in 1972 \cite{zeck}.

\bigskip

\begin{theorem}[Zeckendorf's theorem]
Any positive integer $n$ can be expressed uniquely as a sum of Fibonacci numbers 
\begin{equation}
\label{zthm}
n = \sum_{i \geq 0} a_i F_{i+2}
\end{equation}
with $a_i \in \{0,1\}$ and $a_i a_{i+1} = 0$ for all $i$.
\end{theorem}

\bigskip
We identify the integer $n$ with the finite string $(a_0, a_1, \dots )$ given in equation (\ref{zthm}) and call this the Zeckendorf representation of the integer $n$. For example, 
$$
17 = F_2 + F_ 4 + F_7 = (1, 0, 1, 0, 0, 1).
$$

In 1981, Chung and Graham introduced a representation system which only used the even Fibonacci numbers \cite{CG1981}.

\bigskip

\begin{theorem}[Chung and Graham]
\label{cglemma}
Every non-negative integer $n$ can be uniquely represented as a sum
\begin{equation}
\label{cgthm}
n = \sum_{i \geq 0} a_i F_{i+2} \,\,\, \text{ where } a_i \in \{0, 1, 2 \} \text{ and } \,\, a_i = 0 \text{ if $i$ is odd}
\end{equation}
so that, if $ i$ and $j$ are even with $a_i = a_j = 2$ and $i < j$, then there is some even $k: i < k< j$, for which $a_k = 0$.
\end{theorem}

\bigskip

We can again identify the integer $n$ with the finite string $(a_0, a_1, \dots )$ given in equation (\ref{cgthm}).  Following the paper by Chu, Kanji and Vasseur \cite{Chu:2025aa}, we will call this the Chung-Graham decomposition or representation of the integer $n$.

The purpose of this paper is to show the relationship between the Zeckendorf and Chung-Graham representations using the software package {\tt Walnut}. Hamoon Mousavi, who wrote the  {\tt Walnut} program, has provided an introductory article \cite{Mousavi:2016aa}. Further information about {\tt Walnut} can be found in Shallit's book \cite{Shallit:2022}.

The automata in this paper will accept integer representations in least significant digit first (lsd) order unless otherwise stated. 

\bigskip

\section{Chung-Graham representation in  {\tt Walnut}.}

{\tt Walnut} has a facility for implementing number system representations which requires the user to provide {\tt Walnut} with two automata, one to describe the language of the representation and the other to implement addition of two integers in the representation. 

\begin{figure}[htbp]
   \begin{center}
    \includegraphics[width=4in]{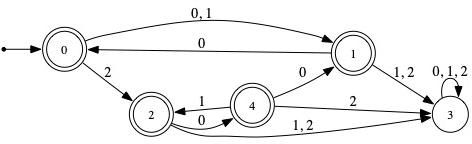}
    \end{center}
    \caption{Automaton which accepts valid Chung-Graham representations.}
    \label{lsdcg}
\end{figure}

The automaton for the language of Chung-Graham representations is shown in Figure \ref{lsdcg}. It takes as input a string of digits from the set $\{0,1,2 \}$ and accepts the string if and only if it is a valid Chung-Graham representation of an integer. We give it the name {\tt cgval}.

\bigskip

\begin{figure}[htbp]
   \begin{center}
    \includegraphics[width=6in]{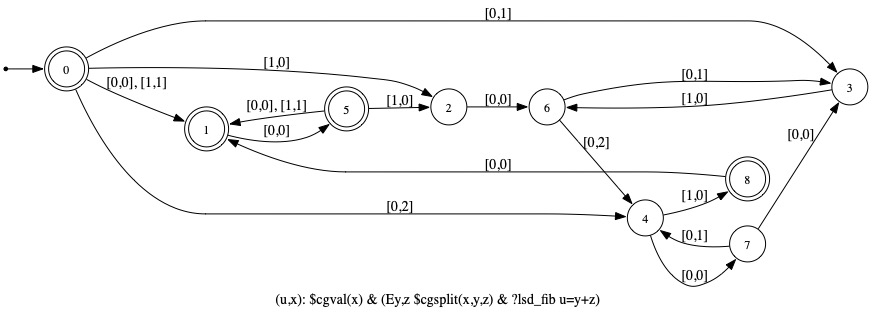}
    \end{center}
    \caption{Automaton which converts between the Zeckendorf and Chung-Graham representations.}
    \label{fibcg}
\end{figure}

\bigskip
Before constructing the automaton which implements addition of Chung-Graham representations, we present an automaton which converts between a valid Chung-Graham representation and the Zeckendorf representation of the same integer. Given a string $u$ defined over the alphabet  $\{0,1,2 \}$, we split $u$ into two strings $v$ and $w$ which are each defined over the alphabet  $\{0,1\}$. This is done using a regular expression which we call {\tt cgsplit}:

\begin{verbatim}
reg cgsplit {0,1,2} {0,1} {0,1} "([0,0,0]|[1,0,1]|[2,1,1])*":
\end{verbatim}
The strings $u, v, w$ satisfy the equation $u=v+w$ when considered as integers via the representations in equations (\ref{zthm}) and (\ref{cgthm}). If $u$ is a valid Chung-Graham representation, then $v$ and $w$ are valid Zeckendorf representations. We can then calculate the sum of $v$ and $w$ using the addition provided for Zeckendorf representations. The result of this addition is the Zeckendorf representation of $u$. The following {\tt Walnut} command performs these actions. The resultant automaton, which we call {\tt fibcg},  is shown in Figure \ref{fibcg}. Its input consists of two integers, which are read in parallel. The first, $u$, is a valid Zeckendorf representation and the second, $x$, is a valid Chung-Graham representation. It accepts the pair if $u = x$ as integers.

\begin{verbatim}
def fibcg "$cgval(x) & (Ey,z $cgsplit(x,y,z) & ?lsd_fib u=y+z)":
\end{verbatim}

\bigskip

We calculated a potential automaton for addition in the Chung-Graham system using experimental data. It is too large to display here, consisting of 33 states. We will call this automaton {\tt cgadd}. Its input consists of three strings $x$, $y$, and $z$ defined over the alphabet $\{0,1,2 \}$. It accepts the triple $(x,y,z)$ if and only if $x$, $y$ and $z$ are valid Chung-Graham representations and $x + y = z$ as integers. 

Ideally, for every $x$ and $y$ there is a unique $z$ such that $x+y=z$. In order for {\tt cgadd} to find a suitable $z$ we need to pad the end of $x$ and $y$ with $0$'s. This ensures that a suitable $z$ exists that is the same length as $x$ and $y$. The regular expression {\tt cg0} defines the set of strings which end in two $0$'s at the most significant digit end:
\begin{verbatim}
reg cg0 {0,1,2} "[0|1|2]*00":
\end{verbatim}

We also need a way of determining when two strings are the same. This is done by the regular expression {\tt cgeq} below:
\begin{verbatim}
reg cgeq {0,1,2} {0,1,2} "([0,0]||[1,1]|[2,2])*":
\end{verbatim}

We can now show that, when $x$ and $y$ are padded, there is a unique $z$ such that {\tt cgadd} accepts the triple $\{x, y, z \}$. The following {\tt Walnut} commands both return TRUE. The first shows the existence of $z$ and the second shows uniqueness.
\begin{verbatim}
eval tcg "Ax,y $cg0(x) & $cg0(y) & $cgval(x) & $cgval(y) 
    => (Ez $cgadd(x,y,z))":
eval tcg "Aw,x,y,z $cg0(x) & $cg0(y) & $cgval(x) & $cgval(y) 
    & $cgadd(x,y,z) & $cgadd(x,y,w) => $cgeq(z,w)":
\end{verbatim}

Padding $x$ and $y$ with one $0$ is not enough to establish existence of $z$, though it is enough for uniqueness.

We will use {\tt Walnut's} pre-installed addition of Zeckendorf representations to show that {\tt cgadd} performs addition correctly. We do this by comparing the result of adding two integers using {\tt cgadd} with the result obtained by adding the same integers using the Zeckendorf representation of the integers. The following {\tt Walnut} command makes the comparison and returns TRUE. 
\begin{verbatim}
eval tcg "Au,v,w,x,y,z ($fibcg(?lsd_fib u,x) & $fibcg(?lsd_fib v,y) 
    & $fibcg(?lsd_fib w,z)) => (?lsd_fib u+v=w <=> $cgadd(x,y,z))":
\end{verbatim}
 
\bigskip

\section{Converting strings into Chung-Graham form}

In this section we construct an automaton which converts an arbitrary string $s = (a_0, a_1, \dots)$ consisting of digits from $\{ 0,1,2 \}$ into the Chung-Graham representation of the positive integer associated with $s$ via equation (\ref{cgarb}):
\begin{equation}
\label{cgarb}
n = \sum_{i \geq 0} a_i F_{i+2}.
\end{equation}

\begin{figure}[htbp]
   \begin{center}
    \includegraphics[width=4in]{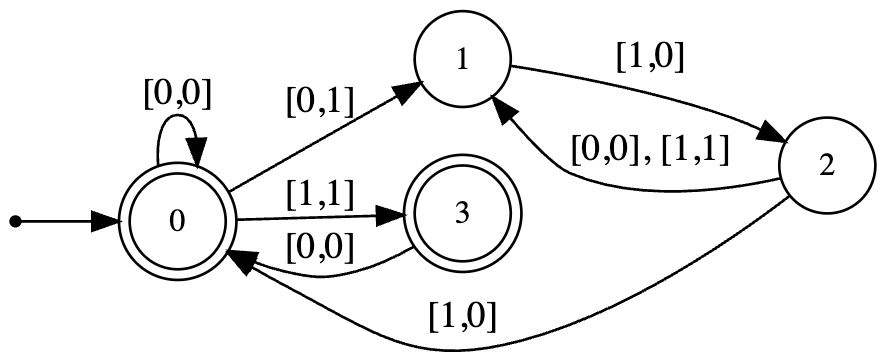}
    \end{center}
    \caption{Automaton which converts an msd binary string into the Zeckendorf representation.}
    \label{bfsc}
\end{figure}

Figure \ref{bfsc} shows an automaton which converts an arbitrary sum of Fibonacci numbers into the Zeckendorf representation of the same integer. The automaton, which is taken from  Figure 1 in the paper \cite{Shallit:2023aa}, takes two binary strings, $x$ and $y$, in parallel and accepts the input if and only if the strings represent the same integer and $y$ is in Zeckendorf format. It takes the input in most significant digit first (msd) order.

\bigskip

\bigskip

\begin{figure}[htbp]
   \begin{center}
    \includegraphics[width=5in]{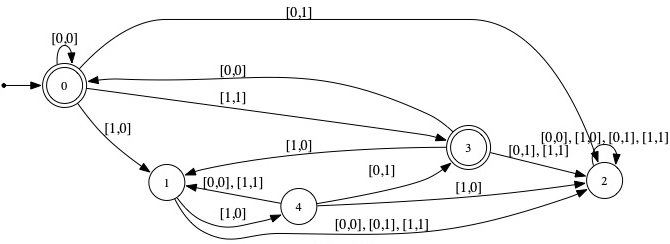}
    \end{center}
    \caption{Automaton which converts an lsd binary string into the Zeckendorf representation.}
    \label{frlsd}
\end{figure}

We can use {\tt Walnut} to convert the automaton at Figure \ref{bfsc}, which we call {\tt fibrepmsd}, into one that uses lsd format. The following command produces an automaton, {\tt fibreplsd}, which takes two binary strings, $x$ and $y$, in parallel lsd format, and accepts the input if and only if the strings represent the same integer and $y$ is in Zeckendorf format.  It is shown at Figure \ref{frlsd}.

\bigskip
\begin{verbatim}
reverse fibreplsd fibrepmsd:
\end{verbatim}

\bigskip

\bigskip

We now consider arbitrary strings, $( a_0, a_1, \dots )$, consisting of the integers $0, 1$, and $2$. Such a string is identified with an integer $n$ via equation (\ref{cgarb}).

We create an automaton which takes as input two such strings, $x$ and $y$, in parallel lsd format, and accepts the input if and only if the strings represent the same integer and $y$ is a valid Chung-Graham representation. The following {\tt Walnut} command splits an arbitrary string $r$ into two strings, $s$ and $t$, which are defined over $\{ 0,1 \}$. It then converts $s$ and $t$ into Zeckendorf forms $u$ and $v$ and takes their sum, producing $w$, which is in Zeckendorf format. Finally, it converts $w$ into Chung-Graham form. The resulting automaton has 42 states.
\begin{verbatim}
def cgrep "Es,t,u,v,w $cgsplit(r,s,t) & $fibreplsd(s,u) & 
    $fibreplsd(t,v) & ?lsd_fib w=u+v & $fibcg(w,z)": 
\end{verbatim}

\bigskip

\section{Synchronisation}
A sequence (or function over $\mathbb{N}$), $f$, is called synchronised if there is an automaton which takes as input a pair of integers $(n, s)$, written in some number system(s), and accepts the pair if and only if $f(n) = s$. If both $n$ and $s$ are written in base $k$, then the function is called $k$-synchronised. If $n$ is written in base-$k$ and $s$ is written in base-$l$, then the function is called $(k, l)$-synchronised. If $n$ and $s$ are written in Zeckendorf representation, then the function is called Fibonacci-synchronised and so on. Since there is an automaton {\tt fibcg} which converts between the the Zeckendorf and Chung-Graham representations, a function is Fibonacci synchronised if and only if it is Chung-Graham-synchronised. We will use {\tt Walnut} to show how to obtain a Chung-Graham-synchronising automaton from a Fibonacci-synchronising automaton. 

\begin{figure}[htbp]
   \begin{center}
    \includegraphics[width=4in]{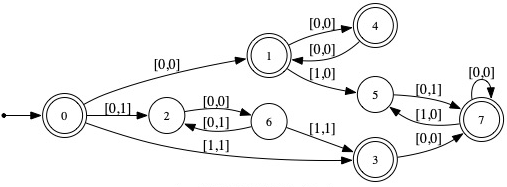}
    \end{center}
    \caption{Synchronised automaton for the function $\lfloor \phi n \rfloor$  in Zeckendorf form.}
    \label{phinlsd}
\end{figure}

The function $f(n) = \lfloor \phi n \rfloor$ is Fibonacci synchronised, as shown in \cite{Shallit:2022} (Theorem 10.11.1). Here $\phi = (1 + \sqrt{5})/2$ is the golden ratio. The automaton in Figure \ref{phinlsd} is an lsd-first synchronising automaton for the function.

\bigskip

If {\tt phinlsd} is the Fibonacci-synchronising automaton for $f$, we can use {\tt Walnut} to construct an automaton which synchronises $f$ when its inputs are written in Chung-Graham form. This is done by the command below, in which the directive {\tt ?lsd\_cg} tells {\tt Walnut} that the inputs are in lsd-Chung-Graham form unless otherwise indicated.

\begin{verbatim}
def cgphin "?lsd_cg Ex,y $phinlsd(?lsd_fib x, ?lsd_fib y) & 
    $fibcg(?lsd_fib x, u) & $fibcg(?lsd_fib y,v)":
\end{verbatim}

\begin{figure}[htbp]
   \begin{center}
    \includegraphics[width=6in]{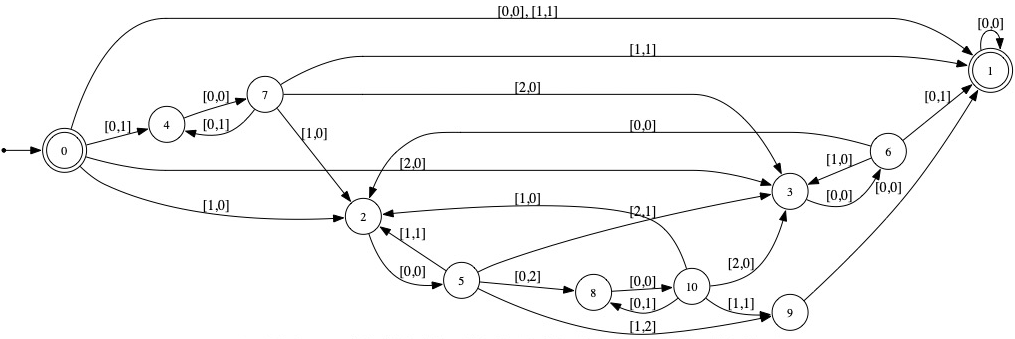}
    \end{center}
    \caption{Synchronised automaton for the function $\lfloor \phi n \rfloor$  in Chung-Graham form.}
    \label{cgphin}
\end{figure}

Figure \ref{cgphin} shows the resulting automaton.

\bigskip

\bibliographystyle{plain}
\begin{small}
\bibliography{CGZ}

\begin{thebibliography}{1}

\bibitem{Chu:2025aa}
Hung~Viet Chu, Aney~Manish Kanji, and Zachary~Louis Vasseur.
\newblock Fixed-term decompositions using even-indexed {F}ibonacci numbers.
\newblock {\em arXiv}, 01 2025.

\bibitem{CG1981}
F.~Chung and R.~L. Graham.
\newblock On irregularities of distribution in finite and infinite sets.
\newblock {\em Colloq. Math. Soc. J{\'a}nos Bolyai}, 37:181--222, 1981.

\bibitem{Gilson:2020aa}
Amelia Gilson, Hadley Killen, Tam{\'a}s Lengyel, Steven~J. Miller, Nadia Razek,
  Joshua~M. Siktar, and Liza Sulkin.
\newblock Zeckendorf's theorem using indices in an arithmetic progression.
\newblock {\em arXiv}, 05 2020.

\bibitem{lekk}
C.~G. Lekkerkerker.
\newblock {V}oorstelling van natuurlijke getallen door een som van {F}ibonacci.
\newblock {\em Simon Stevin}, 29:190--195, 1952.

\bibitem{Mousavi:2016aa}
Hamoon Mousavi.
\newblock Automatic theorem proving in {{\tt Walnut}}.
\newblock {\em arXiv}, 2016.

\bibitem{Shallit:2022}
Jeffrey Shallit.
\newblock {\em The Logical Approach to Automatic Sequences}.
\newblock Cambridge University Press, Sep 2022.

\bibitem{Shallit:2023aa}
Jeffrey Shallit and Sonja~Linghui Shan.
\newblock A general approach to proving properties of {F}ibonacci
  representations via automata theory.
\newblock In {\em Electronic Proceedings in Theoretical Computer Science},
  volume 386, pages 228--24. Open Publishing Association, 2023.

\bibitem{shan2024}
Sonja~Linghui Shan.
\newblock Proving properties of {F}ibonacci representations via automata
  theory.
\newblock Master's thesis, University of Waterloo, 2024.

\bibitem{zeck}
E.~Zeckendorf.
\newblock {R}epr{\'e}sentation des nombres naturels par une somme de nombres de
  {F}ibonacci ou de nombres de {L}ucas.
\newblock {\em Bull. Soc. R. Sci. Li{\`e}ge}, 41:179--182, 1972.

\end{thebibliography}
\end{small}

\end{document}